\documentclass{article}
\usepackage[utf8]{inputenc}

\usepackage{algorithm}
\usepackage{algorithmicx}
\usepackage{algpseudocode}
\usepackage{amsmath}
\usepackage{amssymb}  
\usepackage{amsthm}
\usepackage{arydshln}
\usepackage{bbold}
\usepackage{bm}
\usepackage{booktabs}
\usepackage{bookmark}
\usepackage{calc}              
\usepackage{extarrows}
\usepackage{float}
\usepackage{geometry}
\usepackage{graphicx}          
\usepackage{latexsym}
\usepackage{mathtools}
\usepackage{multirow}          
\usepackage{subcaption}
\usepackage{tikz}              
\usepackage{xcolor}
\usetikzlibrary{shapes.geometric, arrows, chains}

\definecolor{Violet}{named}{violet}

\title{Multicontinuum homogenization  in perforated domains}

\date{}

\usepackage{authblk}

\hypersetup{colorlinks=true,
            linkcolor=blue,
            anchorcolor=blue,
            citecolor=blue}

\begin{document}
\author[1]{Wei Xie}
\author[2]{Yalchin Efendiev}
\author[3]{Yunqing Huang}
\author[4]{Wing Tat Leung}
\author[5]{Yin Yang}
\affil[1]{National Center for Applied Mathematics in Hunan, Xiangtan University, Xiangtan 411105, Hunan, China. \texttt{xiew@smail.xtu.edu.cn}}
\affil[2]{Department of Mathematics, Texas A\&M University, College Station, TX 77843, USA; \texttt{efendiev@math.tamu.edu}}
\affil[3]{School of Mathematics and Computational Science, Xiangtan University, Xiangtan 411105, Hunan, China \texttt{huangyq@xtu.edu.cn}}
\affil[4]{Departament of Mathmatics, City University of Hong Kong, Hong Kong; \texttt{sidnet123@gmail.com}}
\affil[5]{School of Mathematics and Computational Science, Xiangtan University, Xiangtan 411105, Hunan, China \texttt{yangyinxtu@xtu.edu.cn}}

\maketitle

\begin{abstract}

In this paper, we develop a general framework for multicontinuum 
homogenization in perforated domains. 
The simulations of problems in perforated domains are expensive and, 
in many applications, coarse-grid macroscopic models are developed. 
Many previous approaches include homogenization, multiscale finite 
element methods, and so on. In our paper, we design multicontinuum homogenization based on our recently proposed framework. 
In this setting, we distinguish different spatial regions in 
perforations based on their sizes. For example, very thin 
perforations are considered as one continua, while larger 
perforations are considered as another continua. 
By differentiating perforations in this way, we are able to 
predict flows in each of them more accurately. We present a framework
by formulating cell problems for each continuum using 
appropriate constraints for the solution averages and their gradients.
These cell problem solutions are used in a multiscale expansion and 
in deriving novel macroscopic systems for multicontinuum homogenization. 
Our proposed approaches are designed for problems without scale 
separation. We present numerical results for two continuum problems 
and demonstrate the accuracy of the proposed methods.

\end{abstract}

\section{Introduction} \label{sec:introduction}

Problems in perforated domains appear in many applications. 
These include subsurface applications, materials science, membranes, 
filters, and so on. Simulations at the pore scale are very expensive 
and require gridding entailing a very large number of degrees of freedom. 
In many applications, researchers would like to perform simulations 
on a coarse grid and obtain models that do not include perforations. 
These coarse grid models have been a topic of interest for many years. 
In this paper, we propose some novel algorithms for simulations in 
perforated domains.

Some of the first approaches for modeling on a coarse grid include 
homogenization methods, e.g., Darcy's law. In these approaches, 
one assigns an effective property to each representative volume 
based on local simulations. These effective properties are then 
used to form coarse-grid equations. Homogenization techniques are 
specifically designed for problems in perforated domains 
\cite{weinan2003heterogenous, henning2009heterogeneous, hillairet2018homogenization, hornung1997homogenization, yosifian1997some}. 
In these techniques, the solution is expanded using a two-scale ansatz, 
and the terms in the expansion are computed via substitution. 
In many applications, the homogenized equation is derived on 
a coarse grid and does not contain oscillations, while the 
homogenized coefficients are computed based on local cell solutions.

An alternative approach is the use of multiscale methods. 
In these approaches, the solution is sought on a coarse grid 
using multiscale basis functions, which are solutions of local 
problems and are computed locally. 
In \cite{hou1997multiscale, le2014msfem, muljadi2015nonconforming}, 
the authors propose Multiscale Finite ELement Method (MsFEM) approaches, 
where a limited number of basis functions via local solutions are computed, 
and the accuracy and robustness of the approach are demonstrated. 
In \cite{efendiev2013generalized, chung2017online, chung2016mixed, chung2017conservative, chung2018multiscale, chung2016generalized}, 
the authors propose multiscale enrichment and design the Generalized 
Multiscale Finite Element Method (GMsFEM). In these approaches, 
multiscale basis construction is proposed and analyzed. 
By adding additional multiscale basis functions on the pore scale, 
the accuracy of the coarse-grid simulation improves. 
In further generalization \cite{chung2018constraint, chung2021convergence}, 
the authors propose the Constraint Energy Minimizing Generalized Multiscale
Finite Element Method (CEM-GMsFEM) approach. In this approach, 
multiscale basis functions are constructed in oversampled regions 
using constraints. It can be shown that the accuracy of these 
approaches is independent of small scales.

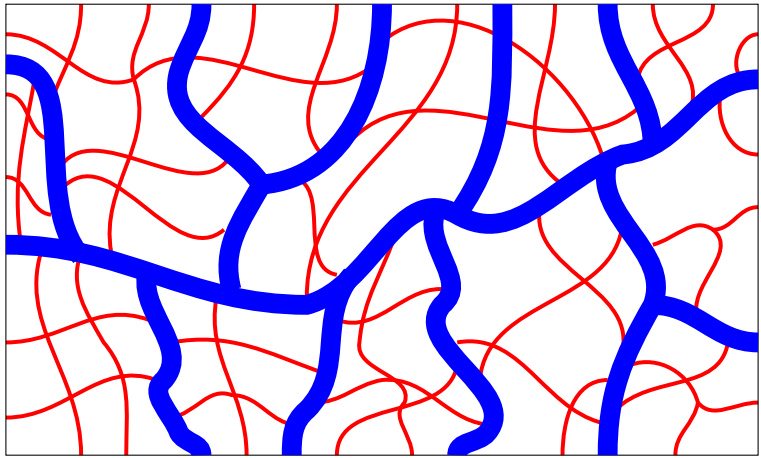
\begin{figure}[H]
\centering
\begin{tikzpicture}
\draw[color=red, line width=1.5pt] (0,3.7) to[out=0,in=175] (0.6,3.2);
\draw[color=red, line width=1.5pt] (0,4.8) to[out=0,in=175] (0.6,4.2);
\draw[color=red, line width=1.5pt] (0.2,2.8) to[out=98,in=260] (0.4,5.1);
\draw[color=red, line width=1.5pt] (0.5,2.8) to[out=250,in=90] (1,0);
\draw[color=red, line width=1.5pt] (0,0.5) to[out=0,in=160] (2,1);
\draw[color=red, line width=1.5pt] (0,1.5) to[out=0,in=160] (1.9,1.8);
\draw[color=red, line width=1.5pt] (0,5.6) to[out=0,in=200] (1.7,5) to[out=45,in=220] (4.9,5.2);
\draw[color=red, line width=1.5pt] (1,2.7) to[out=250,in=120] (1.3,1.5) to[out=310,in=90] (1.6,0);
\draw[color=red, line width=1.5pt] (1.4,2.7) to[out=100,in=290] (1.7,5) to[out=105,in=270] (1.9,6);
\draw[color=red, line width=1.5pt] (0.7,3.8) to[out=60,in=240] (3,4);
\draw[color=red, line width=1.5pt] (0.8,3.6) to[out=40,in=220] (2.9,3);
\draw[color=red, line width=1.5pt] (0.6,4.7) to[out=70,in=270] (1,6);
\draw[color=red, line width=1.5pt] (2.2,0.5) to[out=10,in=160] (3,0.8) to[out=0,in=180] (3.9,0.4);
\draw[color=red, line width=1.5pt] (2.2,1.5) to[out=10,in=160] (4.2,1.1);
\draw[color=red, line width=1.5pt] (3.2,0) to[out=90,in=260] (2.8,2.1);
\draw[color=red, line width=1.5pt] (4.1,0.7) to[out=0,in=190] (5,1) to[out=5,in=180] (6.5,0.4);
\draw[color=red, line width=1.5pt] (4.2,4) to[out=130,in=270] (4.4,6);
\draw[color=red, line width=1.5pt] (5.2,3) to[out=252,in=85] (4.7,1.5)
to[out=260,in=70] (5.3,0.7) to[out=225,in=90] (5.4,0);
\draw[color=red, line width=1.5pt] (4.4,1.8) to[out=340,in=170] (5.9,2.2);
\draw[color=red, line width=1.5pt] (4.4,0) to[out=90,in=260] (6,1.2);
\draw[color=red, line width=1.5pt] (4.6,4.2) to[out=45,in=225] (8.5,4.7);
\draw[color=red, line width=1.5pt] (4,2.1) to[out=90,in=270] (6,6);
\draw[color=red, line width=1.5pt] (3.8,3.7) to[out=350,in=170] (4.4,2.4);
\draw[color=red, line width=1.5pt] (2.6,4.5) to[out=45,in=270] (3.3,6);
\draw[color=red, line width=1.5pt] (7.4,0) to[out=90,in=180] (8,1);
\draw[color=red, line width=1.5pt] (7.4,3.6) to[out=145,in=270] (7.3,6);
\draw[color=red, line width=1.5pt] (5,5.5) to[out=40,in=100] (8,3.7);
\draw[color=red, line width=1.5pt] (7.1,3.3) to[out=260,in=80] (8.2,1.4);
\draw[color=red, line width=1.5pt] (6,1.5) to[out=10,in=175] (8.1,0.4);
\draw[color=red, line width=1.5pt] (6.3,1) to[out=85,in=265] (8.2,3.1);
\draw[color=red, line width=1.5pt] (8.3,1.2) to[out=20,in=105] (9.2,0.7) to[out=260,in=90] (8.9,0);
\draw[color=red, line width=1.5pt] (9.2,0.7) to[out=310,in=180] (10,0.7);
\draw[color=red, line width=1.5pt] (9.6,1.5) to[out=260,in=45] (8.1,0.6);
\draw[color=red, line width=1.5pt] (9.5,4.8) to[out=260,in=180] (10,4);
\draw[color=red, line width=1.5pt] (8.6,2.8) to[out=15,in=145] (9.4,3) to[out=0,in=180] (10,3.3);
\draw[color=red, line width=1.5pt] (9.4,3) to[out=320,in=110] (9.2,1.8);
\draw[color=red, line width=1.5pt] (8.4,5) to[out=82,in=270] (9.4,6);
\draw[color=red, line width=1.5pt] (9.8,5) to[out=143,in=180] (10,5.6);
\draw[color=red, line width=1.5pt] (8.6,6) to[out=270,in=70] (9.1,4.6);
\draw (0,0) rectangle (10,6);
\draw[color=blue,line width=7.5pt] (0,2.8) to[out=0, in=180]  (4,2) to[out=15,in=150] (6,3.2)
to[out=330,in=200] (8.2,4) to[out=5,in=180] (10,5);
\draw[color=blue,line width=7.5pt] (0,5.2) to[out=0,in=130] (1,2.64);
\draw[color=blue,line width=7.5pt] (2.6,0) to[out=90,in=300] (2.3,0.3) to[out=110,in=220] (2.1,1)
to[out=50,in=250] (1.9,2.5);
\draw[color=blue,line width=7.5pt] (3,2.2) to[out=105,in=240] (3.4,3.6) to[out=5,in=270] (5,6);
\draw[color=blue,line width=7.5pt] (3.4,3.6) to[out=125,in=255] (2.3,5.1) to[out=73,in=270] (2.6,6);
\draw[color=blue,line width=7.5pt] (6,3.2) to[out=55,in=270] (6.6,6);
\draw[color=blue,line width=7.5pt] (3.8,0) to[out=90,in=220] (4,0.6) to[out=45,in=230] (4.6,2.4);
\draw[color=blue,line width=7.5pt] (6,0) to[out=90,in=230] (6.4,0.3) to[out=50,in=230] (5.8,2)
to[out=40,in=255] (5.7,3.2);
\draw[color=blue,line width=7.5pt] (8,3.9) to[out=255,in=70] (8.6,2) to[out=240,in=90] (8,0);
\draw[color=blue,line width=7.5pt] (8.6,2) to[out=0,in=180] (10,1.5);
\draw[color=blue,line width=7.5pt] (8.6,4.1) to[out=90,in=270] (8,6);
\end{tikzpicture}
\caption{Illustration of two continua (red and blue).}
\label{fig:river}
\end{figure}

In our paper, we utilize the concept of multicontinuum homogenization 
proposed in \cite{efendiev2023multicontinuum}. The main idea of 
multicontinuum homogenization is to formulate a coarse-grid equation 
using constraint cell problems. These approaches borrow some ideas 
from CEM-GMsFEM; however, their goal is to derive macroscopic equations
and identify macroscopic variables that smoothly vary over the spatial region. 
In multicontinuum homogenization, one deals with problems where multiple 
macroscopic media are present. In this case, we consider perforated 
domains with vastly different perforation sizes (see \autoref{fig:river}). 
Such perforated media arise in many applications, where one deals with 
vastly different sizes of regions, for example, blood vessels, 
perforation sizes, fractures, vugs, and so on. In these problems, 
it is more advantageous to separate the two media, as the solution in 
each type of perforation (small and large) can behave drastically 
differently. Indeed, if the effects of these perforations are lumped 
into one average across both small and large perforations, 
then the effects of small perforations will be ignored, 
as their effects are much weaker compared to those in large ones. 
For this reason, we propose a multicontinuum approach for such 
homogenization that can handle problems without scale separation.

Some of the main ingredients of the proposed approach are the use of 
constraint cell problems and multicontinuum homogenization expansion. 
In our approach, we propose cell problem formulations that are 
constrained in different parts of the perforations. More precisely, 
our cell problems constrain averages and gradients of the solution 
in subregions of perforations. These cell problems are used in a 
multiscale multicontinuum expansion. By substituting this expansion 
into the macroscale equation, we derive a system of equations that 
describe the coarse-grid system.

One of the main differences between our approach and multiscale 
methods is that the proposed methods provide a coarse-grid model 
in the form of differential equations. This is because we seek 
smooth functions that can approximate the coefficients in multiscale
numerical approaches and can formulate macroscopic models for these
global coarse-grid smooth solutions. We perform numerical experiments
for Laplace's equation in perforated domains, even though our approach
can be used for other applications. We choose several types
of perforated domains that include thick and small channels, 
representing two continua. We solve the cell problems and compute the 
effective properties. Effective properties are directional and larger
if more channels are in the corresponding direction. We compare the 
coarse-grid solution to the averaged fine-grid solutions. The averages 
are taken in each subregion. Our numerical results show that the errors
are small, and the macroscopic model accurately predicts the averages
of the fine-grid solution.

Our contributions in this paper are as follows:
\begin{itemize}
\item Development of a framework for multicontinuum homogenization 
in perforated regions by identifying various continua regions.
\item Formulation of multicontinuum constraint cell problems and multiscale expansion.
\item Derivation of macroscopic equations in the form of coupled convection-diffusion-reaction equations.

\item Presentation of numerical results for various types of multicontinua media.
\end{itemize}

The paper is organized as follows. 
In Section \ref{sec:preliminaries}, we present preliminaries.
Section \ref{sec:multicontinuumhomogenization} is devoted to the 
description of multicontinua approach. 
The numerical results are presented in Section \ref{sec:numericalresults}.
We make some conclusions in Section \ref{sec:conclusions}.

\section{Preliminaries} \label{sec:preliminaries}
In this section, we present the model problem and 
review some previous work. 
Consider the following equations in a perforated domain,

\begin{equation}
\mathcal{L}(u) = f ~\text{in}~ \Omega^{\epsilon},
\label{eq:pde_perforated}
\end{equation}
subject to some boundary and initial conditions. 
Here, $\epsilon$ denotes multiscale quantities, 
such as domains or variables, and 
$\Omega \subset \mathbb{R}^d$ (with $d=2,3$) 
is a bounded domain. Let $\mathcal{B}^{\epsilon}$ 
be the set of perforations within $\Omega$, 
and define $\Omega^{\epsilon} = \Omega \setminus 
\mathcal{B}^{\epsilon}$. 

To simplify the notations, let $V(\Omega^{\epsilon})$ be the 
appropriate solution space, and 

$$
V_0(\Omega^{\epsilon}) = \{ v \in V(\Omega^{\epsilon}), v = 0 ~\textnormal{on}~\partial \Omega^{\epsilon} \}.
$$
The variational formulation of \eqref{eq:pde_perforated} is 
to find $u \in V(\Omega^{\epsilon})$ such that 

$$
( \mathcal{L}(u), v )_{\Omega^{\epsilon}} = 
( f, v )_{\Omega^{\epsilon}}, \quad
\forall v \in V_0(\Omega^{\epsilon}),
$$
where $( \cdot, \cdot )$ denotes a specific inner product for scalar functions or vector functions in the perforated domain $\Omega^{\epsilon}$.
In the following, we provide two examples for the abstract notations.

\begin{enumerate}
\item For the Laplace operator,
\begin{equation}
\mathcal{L}(u) = -\Delta u,
\label{eq:pde_laplace}
\end{equation}
we assume the boundary conditions is homogeneous Dirichlet. 

\item For Stokes equations, we have
\begin{equation}
\mathcal{L}(u,p) = 
\begin{pmatrix}
\nabla p - \mu \Delta u \\ \nabla \cdot u
\end{pmatrix},
\label{eq:pde_stokes}
\end{equation}
where $\mu$ is the viscosity, $p$ is the fluid pressure, $u$ represents the velocity.
\end{enumerate}

Next, we review some related techniques for multiscale modeling 
for problems in perforated domains. These include
homogenization, MsFEM, GMsFEM techniques and their variations.

\begin{figure}[htbp]
\centering
\begin{tikzpicture}
\draw[line width=1pt] (0,0) grid (4,4);
\draw[fill=Violet] (0.5,0.5) circle (0.1);
\draw[fill=Violet] (0.5,1.5) circle (0.1);
\draw[fill=Violet] (0.5,2.5) circle (0.1);
\draw[fill=Violet] (0.5,3.5) circle (0.1);
\draw[fill=Violet] (1.5,0.5) circle (0.1);
\draw[fill=Violet] (1.5,1.5) circle (0.1);
\draw[fill=Violet] (1.5,2.5) circle (0.1);
\draw[fill=Violet] (1.5,3.5) circle (0.1);
\draw[fill=Violet] (2.5,0.5) circle (0.1);
\draw[fill=Violet] (2.5,1.5) circle (0.1);
\draw[fill=Violet] (2.5,2.5) circle (0.1);
\draw[fill=Violet] (2.5,3.5) circle (0.1);
\draw[fill=Violet] (3.5,0.5) circle (0.1);
\draw[fill=Violet] (3.5,1.5) circle (0.1);
\draw[fill=Violet] (3.5,2.5) circle (0.1);
\draw[fill=Violet] (3.5,3.5) circle (0.1);
\draw[fill=Violet] (-1,1.5) circle (0.1);
\draw[yellow, line width=2pt] (1,1) rectangle (2,2);
\draw[line width=1pt] (-1.5,1) rectangle (-0.5,2);
\draw[|<->|] (-1.6,1) -- (-1.6,2);
\draw[-latex, red, very thick] (1,1.2) -- (-0.5,1.2);
\draw[-latex, very thick] (-1.3,1.7) -- (-1.3,2.6);
\draw[-latex, very thick] (-1,1.5) -- (-0.9,2.6);
\node at (-1.3,2.7)  {$Y^*$};
\node at (-0.9,2.7)  {$B$};
\node at (-1.7,1.5)  {$\epsilon$};
\draw[line width=1pt] (5,0) grid (9,4);
\draw[yellow, line width=2pt] (6,1) rectangle (7,2);
\draw[blue, line width=2pt] (7,2) rectangle (9,4);
\draw[green, line width=2pt] (5,0) rectangle (8,3);
\draw[fill=Violet] (5.5,0.5) circle (0.15);
\draw[fill=Violet] (5.5,1.5) circle (0.1);
\draw[fill=Violet] (5.5,3.5) circle (0.14);
\draw[fill=Violet] (7,1) circle (0.3);
\draw[fill=Violet] (7.5,2.5) circle (0.2);
\draw[fill=Violet] (8.5,0.5) circle (0.2);
\draw[fill=Violet] (8.5,3) circle (0.25);
\draw[fill=Violet] (6,2.5) ellipse (0.2 and 0.3);
\draw[fill=Violet] (7,3.5) ellipse (0.2 and 0.15);
\draw[blue, line width=2pt] (7,3.35) arc (-95:95:0.2 and 0.15);
\draw[yellow, line width=2pt] (7,1.3) arc (85:185:0.3);
\node at (9.7,1.5)  {$K_i$};
\node at (9.8,0.15)  {$K_{i,1}$};
\node at (9.7,3.5)  {$D_j$};
\draw[red, very thick, -latex] (7,1.5) -- (9.4,1.5);
\draw[red, very thick, -latex] (8,0.15) -- (9.4,0.15);
\draw[red, very thick, -latex] (9,3.5) -- (9.4,3.5);
\end{tikzpicture}
\caption{Left: Periodic perforated domain $\Omega^{\epsilon}$.
Right: Non-periodic perforated domain.}
\label{fig:perforateddomain_pnp}
\end{figure}
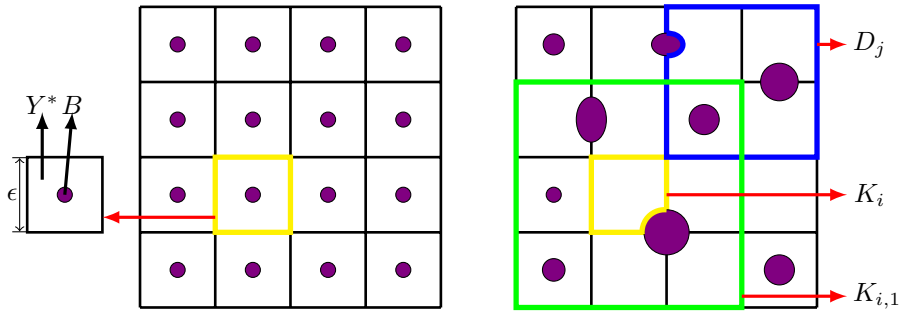

\subsection{Homogenization}

In this subsection, we briefly review homogenization techniques and
refer to \cite{hornung2012homogenization, wolf2022homogenization}. 
We assume the perforations are periodic (see 
\autoref{fig:perforateddomain_pnp}). 
In this context, $\epsilon$ represents the period of the perforations.
For a given cell $Y$, we denote the perforated region as $B$, and 
define the remaining part within $\Omega^{\epsilon}$ 
as $Y^* = Y \setminus B$.

In homogenization, we obtain macroscopic equations using two-scale expansion
that are formulated in the whole domain without perforations
\begin{equation}
\mathcal{L}^* (u^*) = f  \quad \text{in}~ \Omega,
\label{eq:homo_pde}
\end{equation}
where $\mathcal{L}^*$ is the macroscopic operator. The form of the macroscopic operator
depends on microscopic equations. 
For example, for Laplace operator, 
$$
\mathcal{L}^*(u^*) = -\nabla \cdot (a_{ij}^* \nabla u^*),
$$
where
$$
a_{ij}^* = 
\int_{Y^*} \delta_{ij} + \partial_{\boldsymbol{y}_i} w_j(\boldsymbol{y})
~\mathrm{d} \boldsymbol{y},
$$
and $w_j$ satisfied $-\Delta w_j = 0$ in $Y^*$.
On the other hand, for Stokes equations, the macroscopic equations 
differ from the microscopic equations and have the following form
$$
\mathcal{L}^*(u^*, p^*) = 
\begin{pmatrix}
\mathrm{div}(u^*) \\
u^* + a_{ij}^* \nabla p^*
\end{pmatrix}
$$
where $a_{ij}^* = \int_{Y^*} w_j \cdot e^i = 
\int_{Y^*} \nabla w_i : \nabla w_j$, $w_j$ is the cell solution of the following problem,
$$
\begin{cases}
\begin{aligned}
-\mu \Delta w_j + \nabla p_j &= e_j, & \textnormal{in}~ Y^*,\\
\mathrm{div} (w_j) &= 0, & \textnormal{in}~ Y^*, \\
w_j &= 0, & \textnormal{on}~ \partial B.
\end{aligned}
\end{cases}
$$

\subsection{Multiscale finite element method}

Another class of approaches involves the use of multiscale basis 
functions. In these methods, multiscale finite element basis 
functions are constructed to solve the problem on a coarse grid. 
As detailed in \cite{le2014msfem, lu2021uniform}, the MsFEM 
approach can be employed to address non-periodic scenarios, 
as illustrated in \autoref{fig:perforateddomain_pnp}. 
Let $\mathcal{T}^H$ denote a uniform mesh that partitions the 
domain $\Omega$, where $K_i$ represents the coarse cell 
corresponding to the portion within $\Omega^{\epsilon}$, 
i.e., $K_i \subset \Omega^{\epsilon}$. Besides, 
we use $E_j$ to denote the coarse edge, 
and $\mathcal{E}$ is the set of all $E_j$.

For Laplace equations on the fine grid, we need to construct the 
local basis functions as shown in \eqref{eq:msfem_phi}, where 
$\boldsymbol{x}_j$ are coarse vertices in $\Omega^{\epsilon}$.

\begin{equation}
\begin{cases}
\begin{aligned}
-\Delta \phi_j &= 0, &\textnormal{in}~ K_i, \\
\phi_j &= \mu_j, &\textnormal{on}~ \partial K_i.
\end{aligned}
\end{cases}
\label{eq:msfem_phi}
\end{equation}
where $\mu_j$ represents some boundary conditions, as described 
in \cite{hou1997multiscale}. Additionally, more basis functions 
are required to address the absence of coarse vertices not 
included in $\Omega^{\epsilon}$.

\begin{equation}
\begin{cases}
\begin{aligned}
-\Delta \psi_i &= 1, &\textnormal{in}~ K_i, \\
\psi_i &= 0, &\textnormal{on}~ \partial K_i.
\end{aligned}
\end{cases}
\label{eq:msfem_psi}
\end{equation}

\subsection{Generalized multiscale finite element method}

For the Generalized Multiscale Finite Element Method (GMsFEM),
an extension beyond the single basis function paradigm of MsFEM 
is available, as outlined in \cite{chung2016generalized}.
Similar to the MsFEM framework, we employ a rectangular mesh 
$\mathcal{T}^H$ to partition the domain $\Omega$. Here, 
$D_j$ represent the support of multiscale basis functions 
the same as the MsFEM. 

For the snapshot space, we generally have two options. 
The first option is to use all fine-grid basis functions, 
while the other option is to solve the following equations.

\begin{equation}
\begin{cases}
\begin{aligned}
\mathcal{L} (\psi_l^{\textnormal{snap}}) &= 0, 
&\textnormal{in}~ D_i, \\
\psi_l^{\textnormal{snap}} &= \delta_l^h, 
&\textnormal{on}~ \partial D_i.
\end{aligned}
\end{cases}
\end{equation}
The snapshot space in $D_i$ is defined as 
$V^{\textnormal{snap}}(D_i) = \mathrm{span}_l 
\{ \psi_l^{\textnormal{snap}} \}$. 
Then, we need use a spectrum problem to reduce the 
dimension of the local snapshot space, 

\begin{equation}
a_i(\psi_j^i, v) = \lambda_j^i s_i(\psi_j^i, v), \quad 
\forall v \in V^{\textnormal{snap}}(D_i).
\end{equation}
In here, $s_i(u, v) = \int_{D_i} \tilde{\kappa} u v$, 
and $a_i(u, v) = \int_{D_i} \nabla u \cdot \nabla v$ 
for Laplace problem \eqref{eq:pde_laplace}, while 
$a_i(u, v) = \int_{D_i} \nabla u : \nabla v$ 
for Stokes equations \eqref{eq:pde_stokes}.
By selecting the smallest $l_i$ eigenvalues, 
the eigenfunctions will construct the offline space.
At last, the multiscale space is construct by the 
combination of partition unity and offline space,

$$
V_{\textnormal{ms}} = \mathrm{span}_{i,j} \{ \chi_i \psi_j^i \}.
$$

\subsection{Constraint energy minimizing Generalized multiscale finite element method}
For Constraint energy minimizing - Generalized multiscale 
finite element method (CEM-GMsFEM), we need two steps to 
construct the multscale space, refer to 
\cite{chung2021convergence, xie2024cempoisson}.

First, we need to construct a auxiliary space, 
sovling \eqref{eq:aux_cem} in a coarse cell $K_i$,

\begin{equation}
a_i(\phi_j^i, v) = \lambda_j^i s_i(\phi_j^i, v), \quad
\forall v \in V(K_i).
\label{eq:aux_cem}
\end{equation}
and select the first $l_i$ smallest eigenvalues and corresponding eigenfunctions. 
In here, $s_i(u, v) = \int_{K_i} \tilde{\kappa} u v$, and $a_i(u, v) = \int_{K_i} \nabla u \cdot \nabla v$ for Laplace problem \eqref{eq:pde_laplace} and $a_i(u, v) = \int_{K_i} \nabla u : \nabla v$ for Stokes equations \eqref{eq:pde_stokes}.
We can define the local auxiliary space, $V_{\textnormal{aux}}^i = \mathrm{span}\{\phi_i^1, \cdots, \phi_i^{l_i} \}$. 
The global inner product $a(\cdot, \cdot), s(\cdot, \cdot)$ is the sum of each $K_i$, we can also define a global operator $\pi : H_0^1(\Omega^{\epsilon}) \rightarrow V_{\textnormal{aux}}$,

$$
\pi(u) = 
\sum_{i=1}^{N_c} \sum_{j=1}^{l_i} 
\frac{s_i(u, \phi_j^i)}{s_i(\phi_j^i, \phi_j^i)}\phi_j^i,
\quad \forall u \in V.
$$
The multiscale basis functions need to solve a minimization problem.

\begin{enumerate}
\item For Laplace equations,

\begin{equation}
\psi_{j,\textnormal{ms}}^i = \mathrm{argmin}
\{ a(\psi,\psi)+s(\pi(\psi)-\phi_j^i, \pi(\psi)-\phi_j^i)~|~
\psi \in V_0(K_{i,m_i})\}.
\label{eq:r_mini_laplace}
\end{equation}
\item For Stokes problem,

\begin{equation}
\psi_{j,\textnormal{ms}}^i = \mathrm{argmin}
\{ a(\psi,\psi)+s(\pi(\psi)-\phi_j^i, \pi(\psi)-\phi_j^i)~|~
\psi \in V_0(K_{i,m_i}) ~\textnormal{and}~
\nabla \cdot \psi=0 \}.
\label{eq:r_mini_stokes}
\end{equation}
\end{enumerate}

The multiscale space is $V_{\textnormal{ms}} = \mathrm{span}_{i,j} \{\psi_{j,\textnormal{ms}}^i \}.$

\section{Multicontinuum homogenization} \label{sec:multicontinuumhomogenization}

In this section, we introduce the model problem and outline 
the fundamental concept of the multicontinuum homogenization 
method applied to perforated domains.

\begin{figure}[htbp]
\centering
\begin{tikzpicture}
\node[anchor=south west,inner sep=0] at (0,0) {\includegraphics[width=5cm]{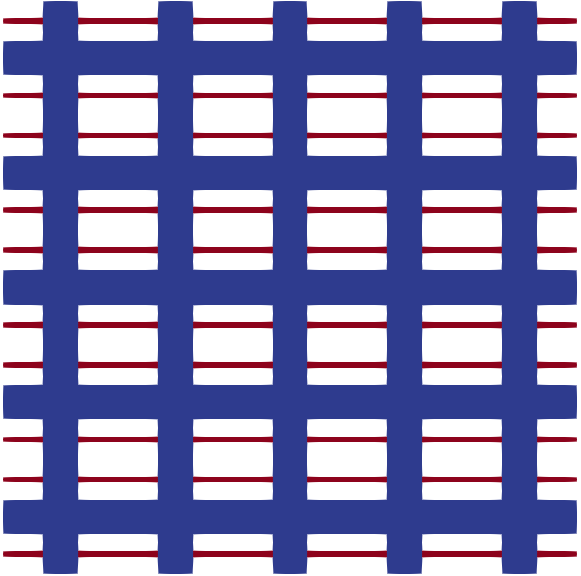}};
\draw[yellow,line width=1.5pt,step=1cm] (0,0) grid (5,5);
\draw[green,line width=2pt] (2,2) rectangle (3,3);
\draw[orange,line width=2pt] (1,1) rectangle (4,4);
\node at (6,2.5)  {$K_{p,1}^{\epsilon}$};
\node at (-1,2.5)  {$K_p^{\epsilon}$};
\draw[line width=1.5pt] (4,1) -- (5.8,2.3);
\draw[line width=1.5pt] (4,4) -- (5.8,2.8);
\draw[line width=1.5pt] (2,2) -- (-1,2.2);
\draw[line width=1.5pt] (2,3) -- (-1,2.8);
\end{tikzpicture}
\caption{Illustration.}
\label{fig:perforated_domain}
\end{figure}
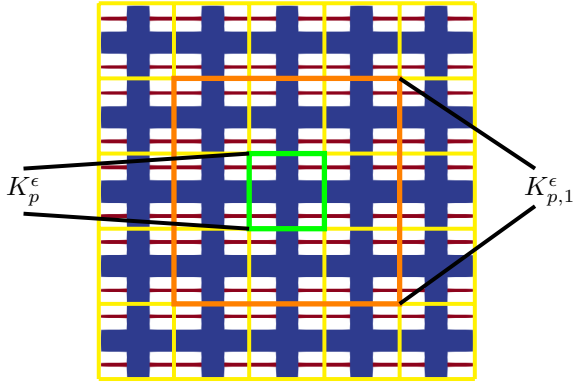

In our study, we define $\mathcal{L}(\cdot) = -\mathrm{div} 
(\kappa \nabla \cdot)$ in Equation \eqref{eq:pde_perforated}, 
subject to homogeneous Dirichlet boundary conditions:

\begin{equation}
\begin{cases}
\begin{aligned}
-\mathrm{div} (\kappa \nabla u) &= f, & \text{ in } \Omega^{\epsilon}, \\
u &= 0, & \text{ on } \partial\Omega^{\epsilon}.
\end{aligned}
\end{cases}
\label{eq:pde}
\end{equation}
Similar to the homogenization example presented in Section 
\ref{sec:preliminaries}, we assume that the perforations 
$\mathcal{B}^{\epsilon}$ exhibit some periodicity in 
Equation \eqref{eq:pde}. 
Here, we use $\epsilon$ to represent the periodicity, 
which corresponds to the width of a single structure.
In applications, channels of different widths may 
possess varying capabilities in transporting flow.
Thus, we categorize continua based on the width of the channels.
In this paper, we focus on the two-continua method, 
as illustrated in Figure \ref{fig:perforated_domain}.
Specifically, we use two distinct colors to denote different 
channels, where blue represents thick channels and red 
indicates thin channels. 
It's worth noting that this method can be easily extended to 
handle multi-channel scenarios.

The weak formulation of Equation \eqref{eq:pde} is given by:

\begin{equation}
a(u, v) = (f, v), \quad \forall v \in H_0^1(\Omega^{\epsilon}),
\label{eq:weak_pde}
\end{equation}
where

$$
a(u,v)=\int_{\Omega^{\epsilon}} \kappa\nabla u\nabla v,\quad
(f,v) =\int_{\Omega^{\epsilon}} fv.
$$

Before presenting the computation, 
let's review the computational mesh. 
We denote a coarse block in the partition of the domain 
$\Omega$ as $K$, with a diameter much larger than 
the heterogeneities, i.e., larger than the smallest 
size of the regions. 
Use $N_c$ to denote the number of coarse grid.
Similar to CEM-GMsFEM, 
we extend the coarse block $K_p$ by $l$ coarse blocks, 
denoted by $K_{p,l}$. We use $K^{\epsilon}$ to 
represent the computational domain within $K$, 
defined as $K^{\epsilon} = \Omega^{\epsilon} \cap K$. 
Additionally, we introduce $K_{p,l}^{\epsilon}$ as the remaining 
part of the oversampled domain $K_{p,l}$ in 
$\Omega^{\epsilon}$, which may contain more local geometric 
information, as illustrated in Figure \ref{fig:perforated_domain}.
To distinguish between different continua, we define the 
characteristic function for continuum $i$ as $\psi_i$, 
satisfying $\psi_i=\delta_{ij}$ within continuum $j$.

In the multicontinuum homogenization method, 
we postulate that the solution $u$ can be expressed as a 
series expansion of macroscopic variables $U_i$. 
Typically, due to the non-periodic nature of geometric 
configurations or coefficients within the operator, 
such as permeability $\kappa$, local information $\phi$ 
must be computed within each coarse block $K_p^{\epsilon}$. 
It is crucial to emphasize that while the variable $U_i$ 
is defined across the global domain $\Omega$, the local 
information $\phi_i$ is confined to $\Omega^{\epsilon}$.
We represent the solution $u$ as follows:

\begin{equation}
u = \phi_i U_i + \phi_i^m \partial_m U_i + \phi_i^{mn} \partial_{mn}^2 U_i + \cdots,
\end{equation}
where $\partial_m = \frac{\partial}{\partial x_m}, 
\partial_{mn}^2 = \frac{\partial^2}{\partial x_m \partial x_n}$.
In general, we consider only the average and the gradient,

\begin{equation}
u \approx \phi_i U_i + \phi_i^m \partial_m U_i.
\end{equation}

\begin{figure}[htbp]
\centering
\begin{tikzpicture}
\node[anchor=south west,inner sep=0] at (0,0) {\includegraphics[width=6cm]{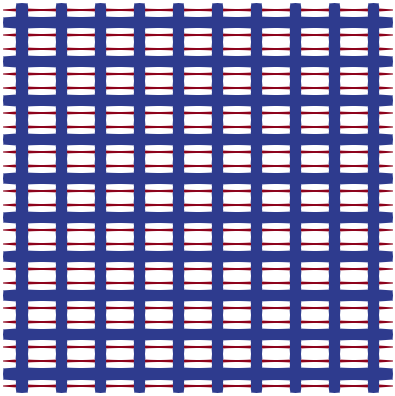}};
\node[anchor=south west,inner sep=0] at (7,3.3) {\includegraphics[width=2.7cm]{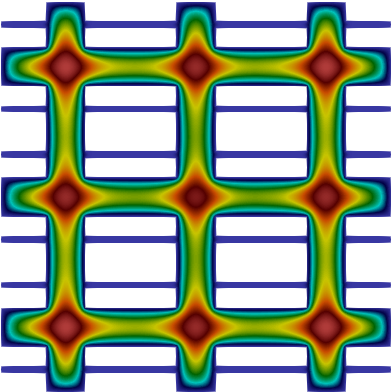}};
\node[anchor=south west,inner sep=0] at (10.5,3.75) {\includegraphics[width=1.5cm]{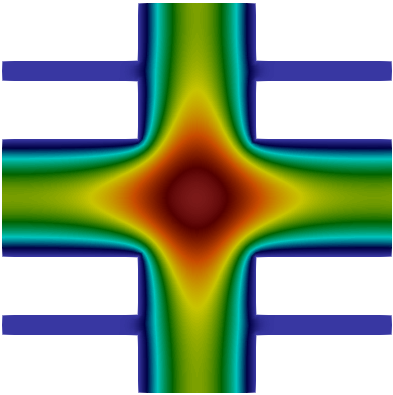}};
\node[anchor=south west,inner sep=0] at (7,0) {\includegraphics[width=2.7cm]{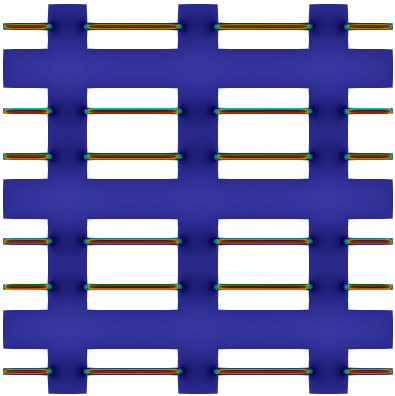}};
\node[anchor=south west,inner sep=0] at (10.5,0.75) {\includegraphics[width=1.5cm]{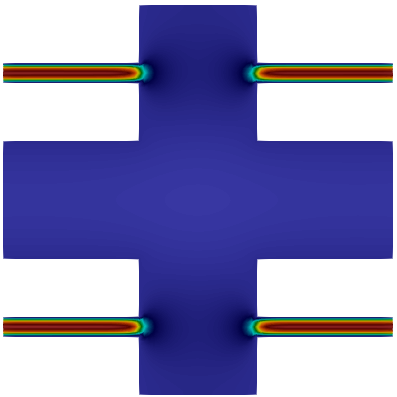}};
\draw[yellow, line width=1.5pt, step=0.6cm] (0,0) grid (6,6);
\draw[orange,line width=1.5pt] (2.4,2.4) rectangle (3,3);
\draw[green,line width=1.5pt] (1.8,1.8) rectangle (3.6,3.6);
\draw[yellow, line width=1.5pt, step=0.9cm, shift={(7,3.3)}] (0,0) grid (2.7,2.7);
\draw[yellow, line width=1.5pt, step=0.9cm, shift={(7,0)}] (0,0) grid (2.7,2.7);
\draw[green,line width=1.5pt] (7,0) rectangle (9.7,2.7);
\draw[green,line width=1.5pt] (7,3.3) rectangle (9.7,6);
\draw[orange,line width=1.5pt] (7.9,0.9) rectangle (8.8,1.8);
\draw[orange,line width=1.5pt] (7.9,4.2) rectangle (8.8,5.1);
\draw[red, very thick, -latex] (3.6,2.3) -- (6.8,1.35);
\draw[red, very thick, -latex] (8.8,1.35) -- (10.5,1.35);
\draw[red, very thick, -latex] (3.6,3.2) -- (6.8,4.65);
\draw[red, very thick, -latex] (8.8,4.65) -- (10.5,4.65);
\end{tikzpicture}
\caption{Illustration of the computational process for a cell problem.}
\label{fig:phi_construction}
\end{figure}
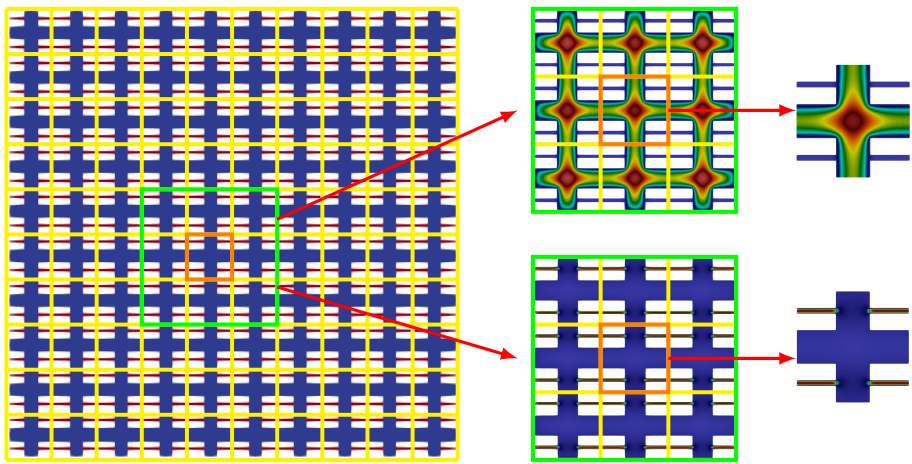

The construction of local information $\phi_i$ is divided into 
$N_c$ local problems, as depicted in 
\autoref{fig:phi_construction}. Here, we denote 
$\phi_i |_{K_p^{\epsilon}} = \phi_i^p$. Subsequently, 
the superscript $p$ will be omitted, as local computation will 
be performed for each coarse block $K_p^{\epsilon}$, i.e., 
$\phi_i |_{K_p^{\epsilon}} = \phi_i$. From previous work 
\cite{efendiev2023multicontinuum, chung2023multicontinuum}, 
it's evident that the cell problem is crucial for obtaining 
an accurate macroscopic equation. These cell problems are 
typically formulated by constraining the original equation 
within an oversampled domain $K_{p,l}^{\epsilon}$. 
Our first cell problem imposes constraints to represent the 
constants in the average behavior of each continuum.

\begin{equation}
\begin{aligned}
\int_{K_{p,l}^{\epsilon}} \kappa \nabla \phi_i \cdot \nabla v - 
\sum_{j,q} \frac{\beta_{ij}^q}{\int_{K_q^{\epsilon}} \psi_j} 
\int_{K_q^{\epsilon}} \psi_j v = 0, \\
\int_{K_q^{\epsilon}} \phi_i \psi_j = \delta_{ij} \int_{K_q^{\epsilon}} \psi_j.
\end{aligned}
\label{eq:phi_avar}
\end{equation}
Our second cell problem imposes constraints to represent 
the linear functions in the average behavior of each continua.

\begin{equation}
\begin{aligned}
\int_{K_{p,l}^{\epsilon}} \kappa \nabla \phi_i^m \cdot \nabla v - 
\sum_{j,q} \frac{\beta_{ij}^{mq}}{\int_{K_q^{\epsilon}} \psi_j} 
\int_{K_q^{\epsilon}} \psi_j v = 0, \\
\int_{K_q^{\epsilon}} \phi_i^m \psi_j = \delta_{ij} \int_{K_q^{\epsilon}} (x_m-c_{mj})\psi_j.
\end{aligned}
\label{eq:phi_grad}
\end{equation}
where $c_{mj}$ satisfy $\int_{K_p^{\epsilon}} (x_m-c_{mj}) \psi_j =0$. 
It's important to note that \eqref{eq:phi_avar} and \eqref{eq:phi_grad}
will use the zero Dirichlet boundary condition, meaning that $\phi$'s 
are all zero on $\partial K_{p,l}^{\epsilon}$, which is consistent with 
the original boundary conditions. The oversampling technique aims to 
remove the boundary effect, but we still need to constrain $\phi$'s 
in $K_p^{\epsilon}$.

Following the above assumption, we only consider the average and 
gradient of macroscopic variables. In particular, we have:

\begin{equation}
\begin{aligned}
u \approx \phi_i U_i + \phi_i^m \partial_m U_i, \\
v \approx \phi_j V_j + \phi_i^m \partial_m V_j.
\end{aligned}
\label{eq:uv_expansion}
\end{equation}

Substituting \eqref{eq:uv_expansion} into \eqref{eq:weak_pde}, 
we obtain the following equation:

\begin{equation}
\begin{aligned}
a(\phi_i U_i, \phi_j V_j) + a(\phi_i U_m, \phi_j^n \partial_n V_j) + 
a(\phi_i^m \partial_m U_i, \phi_j V_j) \\
+ a(\phi_i^m \partial_m U_i, \phi_j^n \partial_n V_j) 
= (f, \phi_j V_j) + (f, \phi_j^n \partial_n V_j).
\end{aligned}
\label{eq:D_aprox}
\end{equation}

Define an local inner product $a_p(u, v) = \int_{K_p^{\epsilon}} 
\kappa \nabla u, \nabla v$.
Considering $U_i$ and $V_i$ as macroscopic variables, 
they can be taken out of the integrals over $K$. Thus, we have:
\begin{equation}
\begin{aligned}
a(\phi_i, \phi_j) U_i V_j + 
a(\phi_i, \phi_j^n) U_i \partial_n V_j 
+a(\phi_i^m, \phi_j) \partial_m U_i V_j  \\
+ a(\phi_i^m, \phi_j^n) \partial_m U_i \partial_n V_j
=(f, \phi_j) V_j +
(f, \partial_n \phi_j) \partial_n V_j.
\end{aligned}
\label{eq:rve_weak}
\end{equation}

In Equation \eqref{eq:rve_weak}, considering that $U_i$ and $V_i$ 
are smooth functions defined in $\Omega$, we derive the following 
macroscopic equation for $U_i$ (in strong form):

\begin{equation}
B_{ji} U_i + B_{ji}^m \partial_m U_i - \partial_n \overline{B}_{ji}^n U_i - \partial_n (B_{ji}^{mn} \partial_m U_i) = b_j,
\label{eq:macro_eq}
\end{equation}
where the coefficients are piecewise-constant vectors or matrices, 
we have:

\begin{equation}
\begin{aligned}
B_{ji} = a(\phi_i, \phi_j), ~
B_{ji}^m = a(\phi_i^m, \phi_j) \\
\overline{B}_{ji}^n = a(\phi_i, \phi_j^n), ~
B_{ji}^{mn} = a(\phi_i^m, \phi_j^n), ~
b_j = (f, \phi_j).
\end{aligned}
\end{equation}

Ultimately, we only need to solve \eqref{eq:macro_eq} 
in $\Omega$ to obtain the macroscopic solution. 
It's important to note:

\begin{enumerate}
\item A suitable oversampling layer can mitigate the boundary effects; 
achieving high accuracy typically requires only one or two layers. 
This is attributed to our knowledge of the interior boundary conditions
and the properties of the macroscopic equation.
\item The concept of the Representative Volume Element (RVE) 
remains applicable in this context.
\end{enumerate}

\section{Numerical results} \label{sec:numericalresults}
In this section, we conduct numerical experiments to examine the 
behavior of two different media with varying conductivity ($\kappa$) 
using four examples. 
The source term will be fixed at $f=5\pi^2 \sin(2\pi x_1) \sin(\pi x_2)$. 
The computation domain is $\Omega=[0,1]^2$. 
Unlike previous studies, we adopt a coarse mesh size $H$ equivalent 
to the periodicity parameter $\epsilon$, which is also the diameter of 
a single structure. 
We will use a fine mesh size of $H/80$ to solve the reference solution.

To measure the efficiency of our method, we define the relative 
$L^2$-error in $\Omega_1$ and the relative $L^2$-error in $\Omega_2$ as:

$$
e_2^{(i)} = \cfrac{\sum_p \left|\frac{1}{|K_p|} \int_{K_p} U_i - \frac{1}{|{K_p^{\epsilon}} \cap \Omega_i|} \int_{{K_p^{\epsilon}}\cap \Omega_i} u \right|^2}{\sum_p \left|\frac{1}{|{K_p^{\epsilon}} \cap \Omega_i|} \int_{{K_p^{\epsilon}}\cap \Omega_i} u\right|^2}.
$$
\begin{figure}[htbp]
\centering
\includegraphics[scale=0.3]{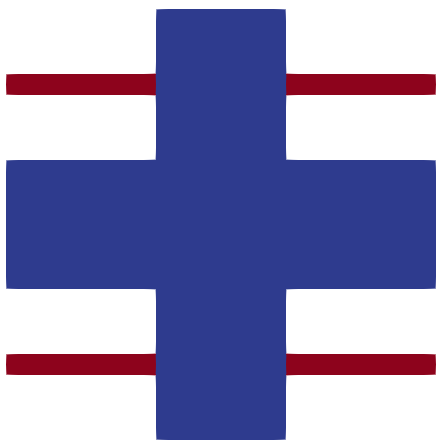}
\includegraphics[scale=0.3]{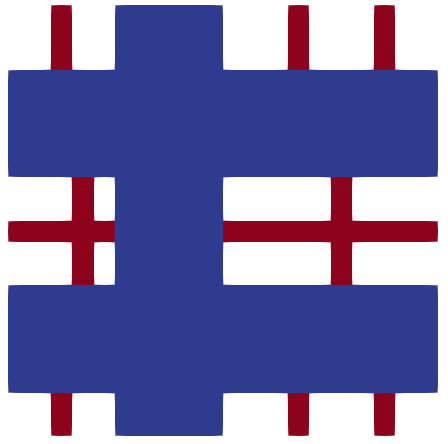}
\caption{Left: single structure 1. Right: single structure 2.}
\label{fig:reference_solution}
\end{figure}

\subsection{Case 1}
\begin{figure}[htbp]
\centering
\includegraphics[scale=0.4]{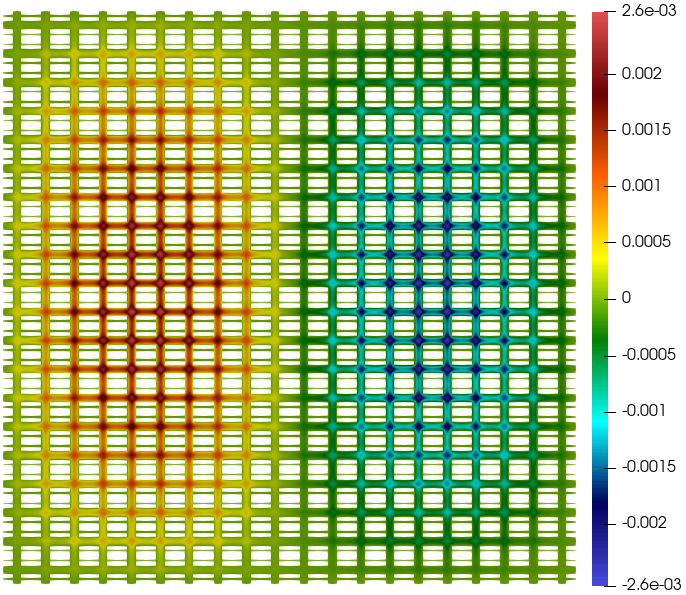}
\caption{Reference solution for Case 1.}
\label{fig:case1_ref}
\end{figure}

\begin{figure}[htbp]
\centering
\includegraphics[scale=0.4]{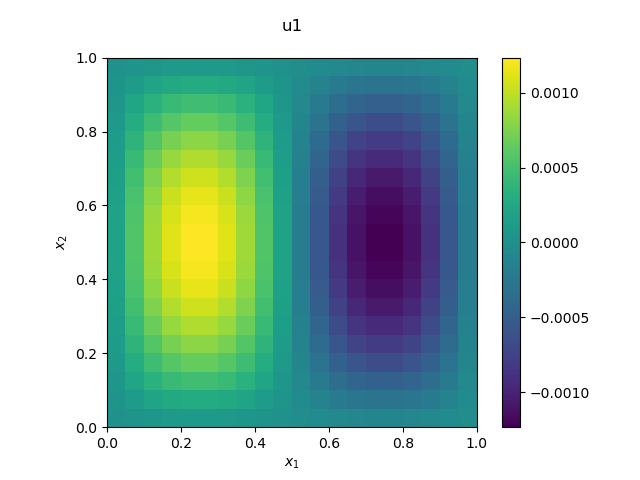}
\includegraphics[scale=0.4]{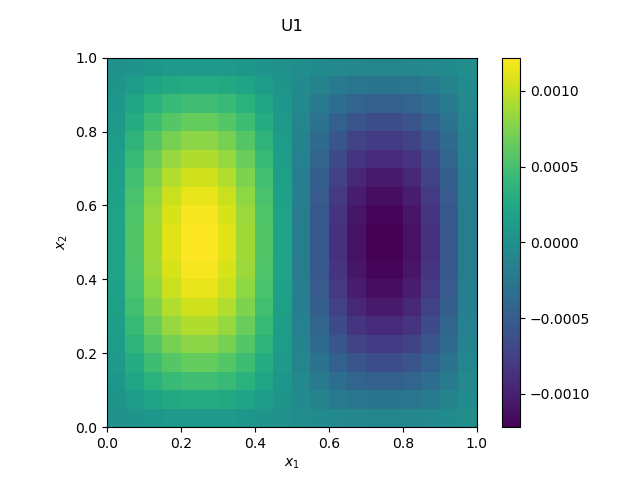}
\includegraphics[scale=0.4]{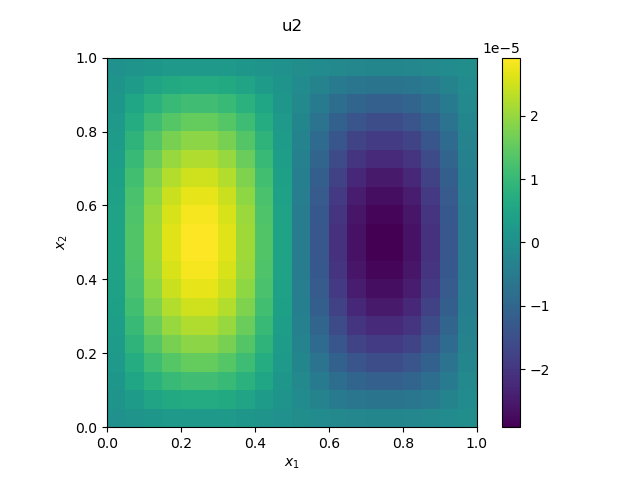}
\includegraphics[scale=0.4]{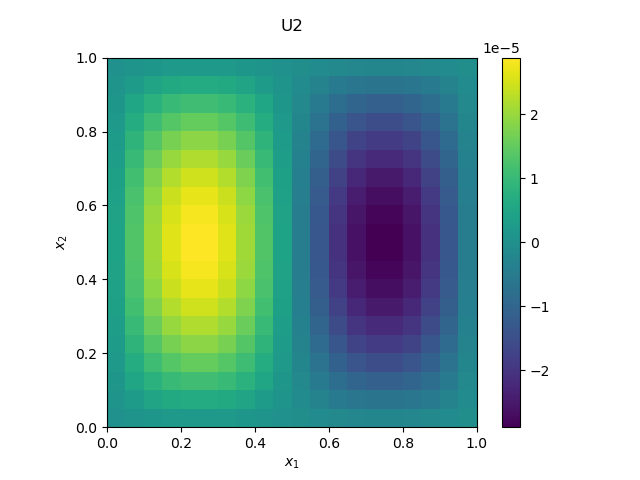}
\caption{Average solution for Case 1. Top Left: Reference averaged solution in $\Omega_1$. Top Right: Multiscale average solution in $\Omega_1$. Bottom Left: Reference averaged solution in $\Omega_2$. Bottom Right: Multiscale average solution in $\Omega_2$.}
\label{fig:case1_uuh}
\end{figure}

In this example, we take $\kappa=1$, and use structure 1 as the single 
structure. We depict the fine-grid solution in \autoref{fig:case1_ref}, 
while the corresponding averaged solution in \autoref{fig:case1_uuh}. 
From \autoref{fig:case1_uuh}, we conclude that our method provides an 
accurate approximation of the averaged solution.
In Table \ref{tab:relative_error_case1}, we observe that the error 
will decay very fast with the coarse mesh size decrease. Besides, we 
find that achieving high accuracy only requires oversampling one coarse 
layer. In Table \ref{tab:case1_Bs}, we present some coefficient in 
macroscopic equation, which can demonstrate the relationship between 
error and oversampling layers.

\begin{table}[htbp]
\centering
\begin{tabular}{|c|c|c|c|c|c|c|}
\hline
\multirow{2}*{$l$} & 
\multicolumn{2}{|c|}{$\epsilon=1/10$} & 
\multicolumn{2}{|c|}{$\epsilon=1/20$} & 
\multicolumn{2}{|c|}{$\epsilon=1/40$}
\\ \cline{2-7}
& $e_2^{(1)}$ & $e_2^{(2)}$ & $e_2^{(1)}$ & $e_2^{(2)}$ & $e_2^{(1)}$ & $e_2^{(2)}$ 
\\ \hline
0 & 4.22e-02 & 1.04e-02 & 3.32e-02 & 5.29e-03 & 3.12e-02 & 4.31e-03
\\ \hline
1 & 2.12e-03 & 1.94e-03 & 1.43e-04 & 1.22e-04 & 9.91e-06 & 7.85e-06
\\ \hline
2 & 2.04e-03 & 1.93e-03 & 1.35e-04 & 1.21e-04 & 8.89e-06 & 7.73e-06
\\ \hline
\end{tabular}
\caption{Relative error in different continuum when use structure 1 and set $\kappa=1$.}
\label{tab:relative_error_case1}
\end{table}

\begin{table}[htbp]
\centering
\begin{tabular}{|c|c|c|c|}
\hline
$l$ & $B_{11}$ & $B_{12}$ & $B_{22}$ \\ \hline
0 & 60.70 & -2.00 & 383.59 \\ \hline
1 & 50.10 & -1.81 & 365.83 \\ \hline
2 & 50.10 & -1.81 & 365.82 \\ \hline
\end{tabular}
\caption{Homogenization coefficient in Case 1.}
\label{tab:case1_Bs}
\end{table}

\subsection{Case 2}

\begin{figure}[htbp]
\centering
\includegraphics[scale=0.4]{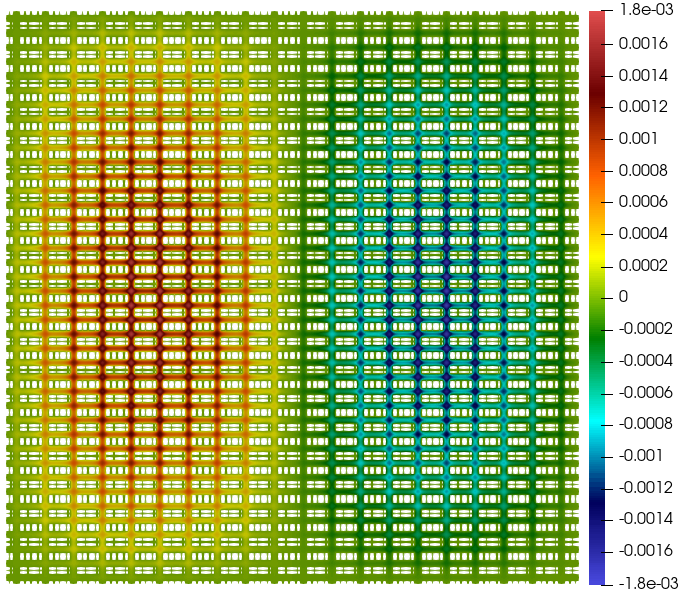}
\caption{Reference solution for Case 2.}
\label{fig:case2_ref}
\end{figure}

\begin{figure}[htbp]
\centering
\includegraphics[scale=0.4]{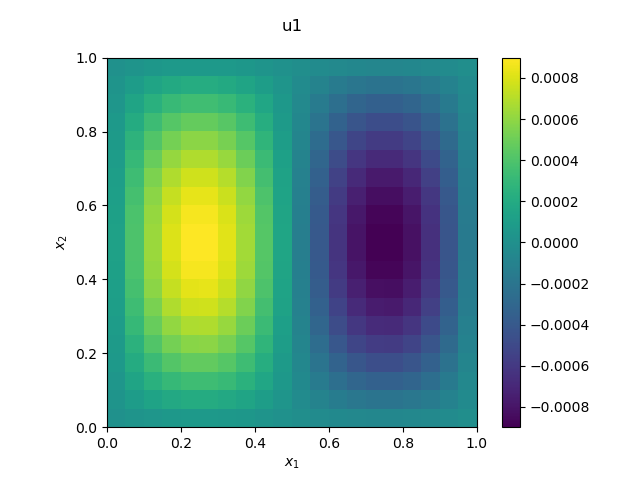}
\includegraphics[scale=0.4]{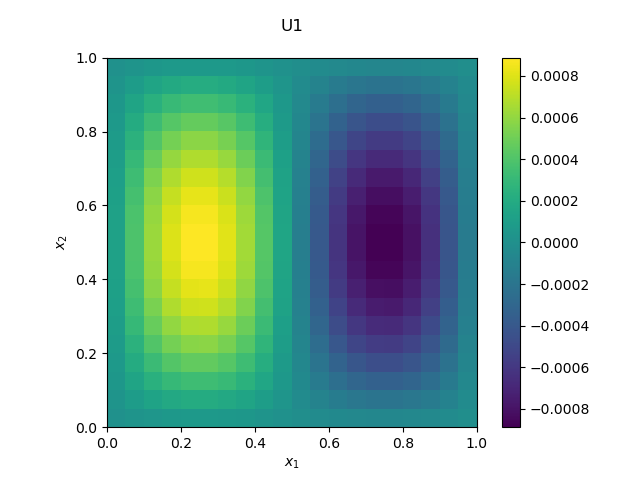}
\includegraphics[scale=0.4]{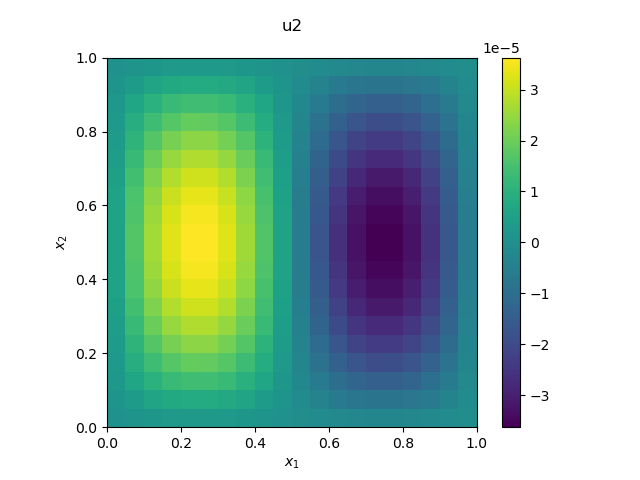}
\includegraphics[scale=0.4]{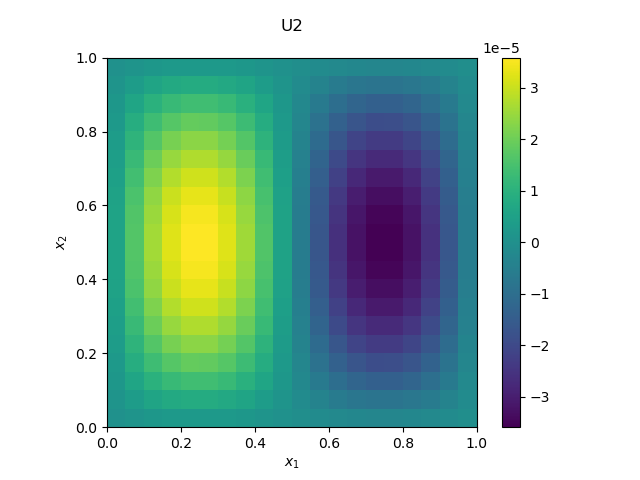}
\caption{Average solution for Case 2. Top Left: Reference averaged solution in $\Omega_1$. Top Right: Multiscale average solution in $\Omega_1$. Bottom Left: Reference averaged solution in $\Omega_2$. Bottom Right: Multiscale average solution in $\Omega_2$.}
\label{fig:case2_uuh}
\end{figure}

In this example, we keep $\kappa=1$ and utilize structure 2 different 
as Case 1. We present the solution obtained from the fine-grid 
approach in \autoref{fig:case2_ref}, alongside the corresponding 
averaged solution displayed in \autoref{fig:case2_uuh}. 
From the analysis of \autoref{fig:case2_uuh}, we infer that our 
method yields a precise approximation of the averaged solution. 
The data in Table \ref{tab:relative_error_case2} indicates a 
rapid decay in error as the coarse mesh size decreases. 
Moreover, we find that achieving high accuracy only requires 
oversampling one coarse layer. Table \ref{tab:case2_Bs} provides 
coefficients in the macroscopic equation, aiding the error will 
converge with only 1 oversampling layer.

\begin{table}[htbp]
\centering
\begin{tabular}{|c|c|c|c|c|c|c|}
\hline
\multirow{2}*{$l$} & 
\multicolumn{2}{|c|}{$\epsilon=1/10$} & 
\multicolumn{2}{|c|}{$\epsilon=1/20$} & 
\multicolumn{2}{|c|}{$\epsilon=1/40$}
\\ \cline{2-7}
& $e_2^{(1)}$ & $e_2^{(2)}$ & $e_2^{(1)}$ & $e_2^{(2)}$ & $e_2^{(1)}$ & $e_2^{(2)}$ 
\\ \hline
0 & 3.81e-02 & 1.17e-02 & 2.93e-02 & 6.36e-03 & 2.73e-02 & 5.30e-03
\\ \hline
1 & 2.26e-03 & 1.94e-03 & 1.57e-04 & 1.23e-04 & 1.13e-05 & 7.97e-06
\\ \hline
2 & 2.17e-03 & 1.92e-03 & 1.47e-04 & 1.21e-04 & 1.02e-05 & 7.77e-06
\\ \hline
\end{tabular}
\caption{Relative error in different continuum when use structure 2 and set $\kappa=1$.}
\label{tab:relative_error_case2}
\end{table}

\begin{table}[htbp]
\centering
\begin{tabular}{|c|c|c|c|}
\hline
$l$ & $B_{11}$ & $B_{12}$ & $B_{22}$ \\ \hline
0 & 101.10 & -6.92 & 494.14 \\ \hline
1 & 84.57 & -5.93 & 462.72 \\ \hline
2 & 84.57 & -5.93 & 462.72 \\ \hline
\end{tabular}
\caption{Homogenization coefficient in Case 2.}
\label{tab:case2_Bs}
\end{table}

\subsection{Case 3}
\begin{figure}[htbp]
\centering
\includegraphics[scale=0.4]{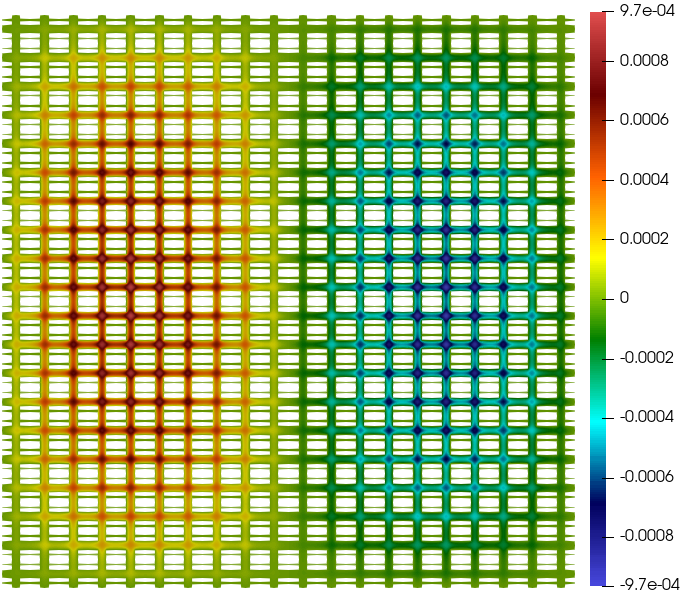}
\caption{Reference solution for Case 3.}
\label{fig:case3_ref}
\end{figure}

\begin{figure}[htbp]
\centering
\includegraphics[scale=0.4]{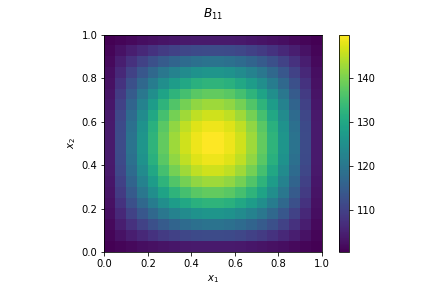}
\includegraphics[scale=0.4]{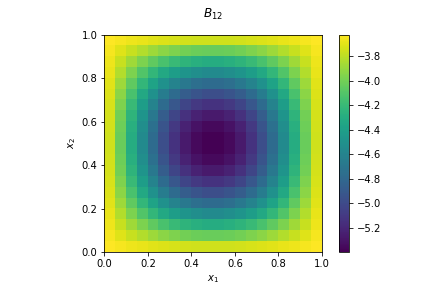}
\includegraphics[scale=0.4]{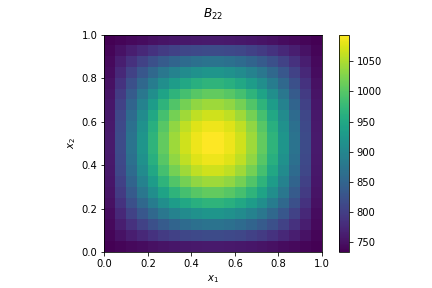}
\caption{Homogenization coefficient for Case 3. From Left to Right and Top to Bottom: $B_{11}, B_{12}, B_{22}$.}
\label{fig:case3_Bs}
\end{figure}

\begin{figure}[htbp]
\centering
\includegraphics[scale=0.4]{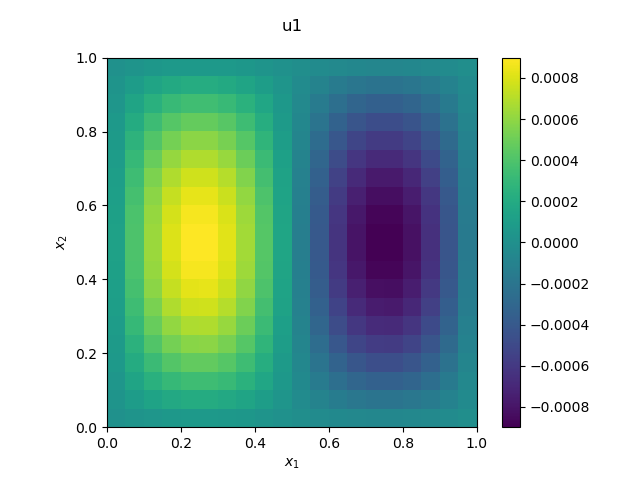}
\includegraphics[scale=0.4]{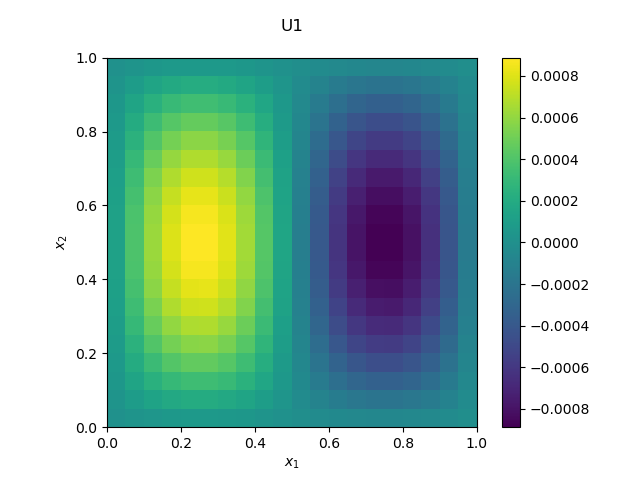}
\includegraphics[scale=0.4]{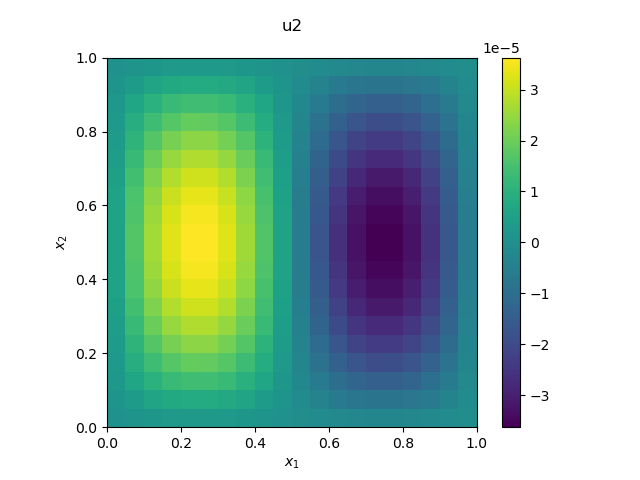}
\includegraphics[scale=0.4]{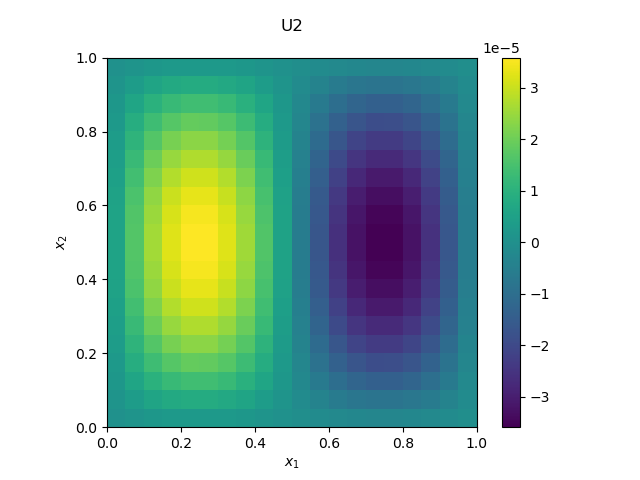}
\caption{Average solution for Case 3. Top Left: Reference averaged solution in $\Omega_1$. Top Right: Multiscale average solution in $\Omega_1$. Bottom Left: Reference averaged solution in $\Omega_2$. Bottom Right: Multiscale average solution in $\Omega_2$.}
\label{fig:case3_uuh}
\end{figure}

In this example, we take a slow variable coefficient, 
$\kappa=2+\sin(\pi x_1)\sin(\pi x_2)$. 
The structure 1 is used for the single structure.
We present the fine-scale solution in \autoref{fig:case3_ref}.
In \autoref{fig:case3_uuh}, we can see that our method can provides 
an accurate approximation of the averaged solution. 
Different with the last two examples, we depicted the 
homogenization coefficient in \autoref{fig:case3_Bs}. 
In Table \ref{tab:relative_error_case3}, we note that the error 
will decay with the reduce the coarse mesh size. 

\begin{table}[htbp]
\centering
\begin{tabular}{|c|c|c|c|c|c|c|}
\hline
\multirow{2}*{$l$} & 
\multicolumn{2}{|c|}{$\epsilon=1/10$} & 
\multicolumn{2}{|c|}{$\epsilon=1/20$} & 
\multicolumn{2}{|c|}{$\epsilon=1/40$}
\\ \cline{2-7}
& $e_2^{(1)}$ & $e_2^{(2)}$ & $e_2^{(1)}$ & $e_2^{(2)}$ & $e_2^{(1)}$ & $e_2^{(2)}$ 
\\ \hline
0 & 4.20e-02 & 1.02e-02 & 3.31e-02 & 5.25e-03 & 3.11e-02 & 4.30e-03
\\ \hline
1 & 2.10e-03 & 1.90e-03 & 1.39e-04 & 1.18e-04 & 9.58e-06 & 7.54e-06
\\ \hline
2 & 2.02e-03 & 1.89e-03 & 1.31e-04 & 1.17e-04 & 8.57e-06 & 7.42e-06
\\ \hline
\end{tabular}
\caption{Relative error in different continuum when use structure 1 and set $\kappa=2+\sin(\pi x_1)\sin(\pi x_2)$.}
\label{tab:relative_error_case3}
\end{table}

\subsection{Case 4}
\begin{figure}[htbp]
\centering
\includegraphics[scale=0.4]{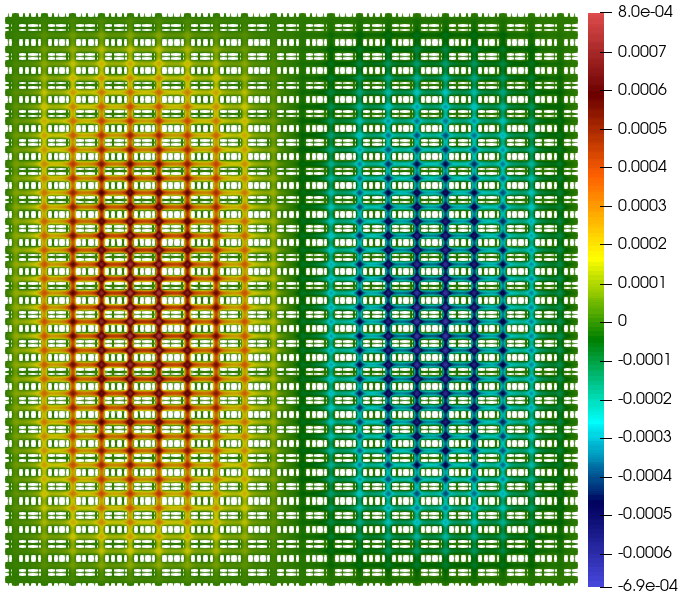}
\caption{Reference solution for Case 4.}
\label{fig:case4_ref}
\end{figure}

\begin{figure}[htbp]
\centering
\includegraphics[scale=0.4]{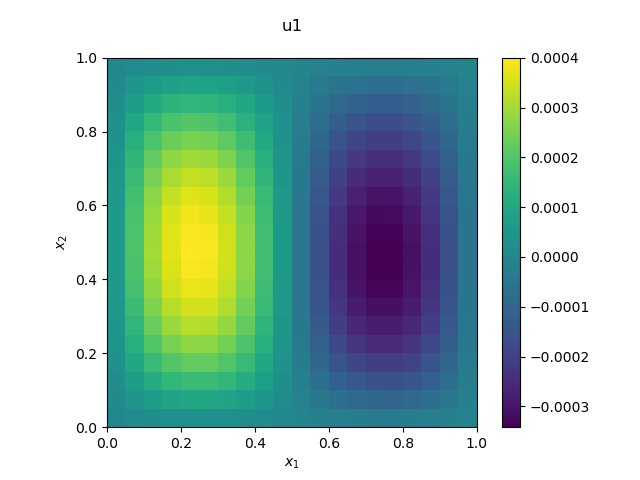}
\includegraphics[scale=0.4]{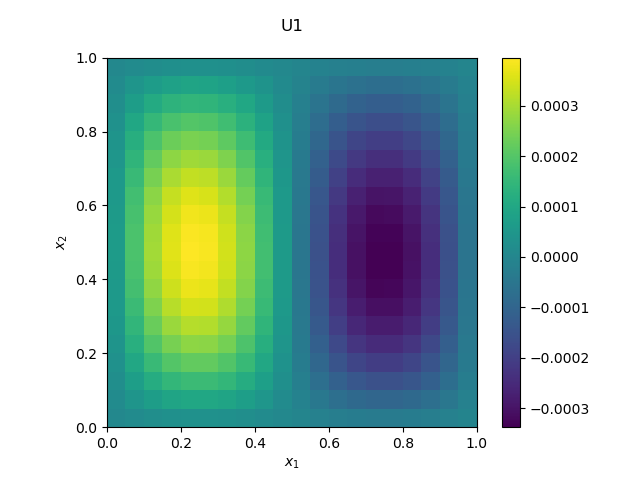}
\includegraphics[scale=0.4]{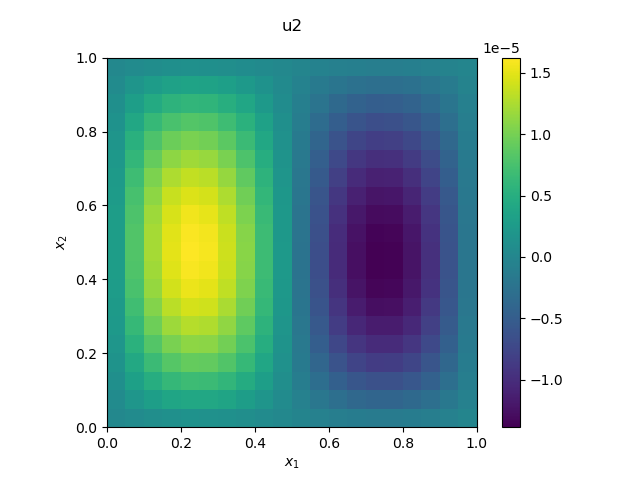}
\includegraphics[scale=0.4]{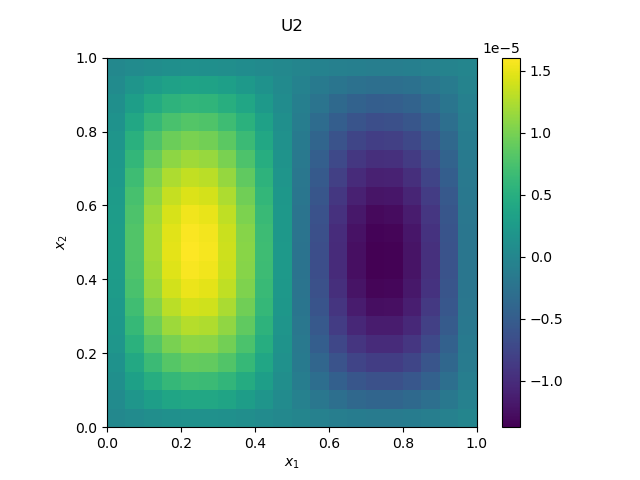}
\caption{Average solution for Case 4. Top Left: Reference averaged solution in $\Omega_1$. Top Right: Multiscale average solution in $\Omega_1$. Bottom Left: Reference averaged solution in $\Omega_2$. Bottom Right: Multiscale average solution in $\Omega_2$.}
\label{fig:case4_uuh}
\end{figure}

\begin{figure}[htbp]
\centering
\includegraphics[scale=0.4]{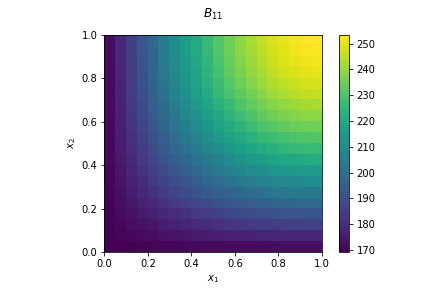}
\includegraphics[scale=0.4]{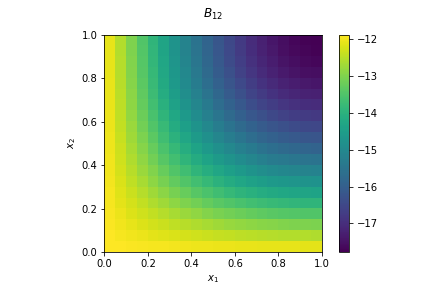}
\includegraphics[scale=0.4]{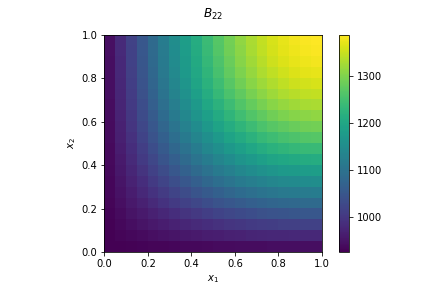}
\caption{Homogenization coefficient for Case 4. From Left to Right and Top to Bottom: $B_{11}, B_{12}, B_{22}$.}
\label{fig:case4_Bs}
\end{figure}

In this example, we set $\kappa=2+\sin(\pi x_1)\sin(\pi x_2)$ 
and utilize Structure 2 as the single structure. 
The fine-grid solution is depicted in \autoref{fig:case4_ref}, 
while the reference averaged solution and multicontinuum 
homogenization solution are presented in \autoref{fig:case4_uuh}. 
The figures further demonstrate the efficiency of our method. 
In Table \ref{tab:relative_error_case4}, we show the error by 
varying the periodic and oversampling layers. We observe that 
the error decays very rapidly when decreasing the periodicity. 
Additionally, in \autoref{fig:case4_Bs}, we show the 
homogenization coefficients.

\begin{table}[htbp]
\centering
\begin{tabular}{|c|c|c|c|c|c|c|}
\hline
\multirow{2}*{$l$} & 
\multicolumn{2}{|c|}{$\epsilon=1/10$} & 
\multicolumn{2}{|c|}{$\epsilon=1/20$} & 
\multicolumn{2}{|c|}{$\epsilon=1/40$}
\\ \cline{2-7}
& $e_2^{(1)}$ & $e_2^{(2)}$ & $e_2^{(1)}$ & $e_2^{(2)}$ & $e_2^{(1)}$ & $e_2^{(2)}$ 
\\ \hline
0 & 3.79e-02 & 1.17e-02 & 2.92e-02 & 6.36e-03 & 2.73e-02 & 5.30e-03
\\ \hline
1 & 2.20e-03 & 1.95e-03 & 1.55e-04 & 1.23e-04 & 1.13e-05 & 7.96e-06
\\ \hline
2 & 2.11e-03 & 1.93e-03 & 1.45e-04 & 1.21e-04 & 1.02e-05 & 7.75e-06
\\ \hline
\end{tabular}
\caption{Relative error in different continuum when use structure 2 and set $\kappa=2+\sin(\frac{\pi x_1}{2})\sin(\frac{\pi x_2}{2})$.}
\label{tab:relative_error_case4}
\end{table}

\section{Conclusions} \label{sec:conclusions}

In this paper, we propose a multicontinuum homogenization approach 
for problems in perforated domains. The perforated regions are 
divided into subregions, where each subregion is treated as a 
separate continua due to their size differences. Typically, 
different continua may have significantly different widths or lengths. 
We formulate constraint cell problems by imposing constraints in 
subregions for the averages of the solutions and their gradients. 
Using the cell solutions, we formulate a homogenization expansion and 
derive macroscopic equations. The resulting macroscopic equations 
consist of a system of equations. We present numerical results by 
considering two continua media with significantly different widths. 
We consider various diffusion scenarios, and our numerical results 
show very good accuracy.

\section*{Acknowledgement}
WX and YY were supported by the National Natural Science Foundation of China Project (12071402, 12261131501), the Project of Scientiﬁc Research Fund of the Hunan Provincial Science and Technology Department (2022RC3022). 
WTL  is supported by 
Early Career Award, Research Grant Council, Project Number: 21307223.
YE would like to thank the partial support from NSF 2208498.

 \bibliographystyle{plain}
\bibliography{references}

\begin{thebibliography}{10}

\bibitem{chung2023multicontinuum}
E~Chung, Y~Efendiev, J~Galvis, and WT~Leung.
\newblock Multicontinuum homogenization. general theory and applications.
\newblock {\em arXiv preprint arXiv:2309.08128}, 2023.

\bibitem{chung2021convergence}
Eric Chung, Jiuhua Hu, and Sai-Mang Pun.
\newblock Convergence of the cem-gmsfem for stokes flows in heterogeneous
  perforated domains.
\newblock {\em Journal of Computational and Applied Mathematics}, 389:113327,
  2021.

\bibitem{chung2018constraint}
Eric~T Chung, Yalchin Efendiev, and Wing~Tat Leung.
\newblock Constraint energy minimizing generalized multiscale finite element
  method.
\newblock {\em Computer Methods in Applied Mechanics and Engineering},
  339:298--319, 2018.

\bibitem{chung2017online}
Eric~T Chung, Yalchin Efendiev, Wing~Tat Leung, Maria Vasilyeva, and Yating
  Wang.
\newblock Online adaptive local multiscale model reduction for heterogeneous
  problems in perforated domains.
\newblock {\em Applicable Analysis}, 96(12):2002--2031, 2017.

\bibitem{chung2016generalized}
Eric~T Chung, Yalchin Efendiev, Guanglian Li, and Maria Vasilyeva.
\newblock Generalized multiscale finite element methods for problems in
  perforated heterogeneous domains.
\newblock {\em Applicable Analysis}, 95(10):2254--2279, 2016.

\bibitem{chung2016mixed}
Eric~T Chung, Wing~Tat Leung, and Maria Vasilyeva.
\newblock Mixed gmsfem for second order elliptic problem in perforated domains.
\newblock {\em Journal of Computational and Applied Mathematics}, 304:84--99,
  2016.

\bibitem{chung2018multiscale}
Eric~T Chung, Wing~Tat Leung, Maria Vasilyeva, and Yating Wang.
\newblock Multiscale model reduction for transport and flow problems in
  perforated domains.
\newblock {\em Journal of Computational and Applied Mathematics}, 330:519--535,
  2018.

\bibitem{chung2017conservative}
Eric~T Chung, Maria Vasilyeva, and Yating Wang.
\newblock A conservative local multiscale model reduction technique for stokes
  flows in heterogeneous perforated domains.
\newblock {\em Journal of Computational and Applied Mathematics}, 321:389--405,
  2017.

\bibitem{efendiev2013generalized}
Yalchin Efendiev, Juan Galvis, and Thomas~Y Hou.
\newblock Generalized multiscale finite element methods (gmsfem).
\newblock {\em Journal of computational physics}, 251:116--135, 2013.

\bibitem{efendiev2023multicontinuum}
Yalchin Efendiev and Wing~Tat Leung.
\newblock Multicontinuum homogenization and its relation to nonlocal
  multicontinuum theories.
\newblock {\em Journal of Computational Physics}, 474:111761, 2023.

\bibitem{henning2009heterogeneous}
Patrick Henning and Mario Ohlberger.
\newblock The heterogeneous multiscale finite element method for elliptic
  homogenization problems in perforated domains.
\newblock {\em Numerische Mathematik}, 113:601--629, 2009.

\bibitem{hillairet2018homogenization}
Matthieu Hillairet.
\newblock On the homogenization of the stokes problem in a perforated domain.
\newblock {\em Archive for Rational Mechanics and Analysis}, 230:1179--1228,
  2018.

\bibitem{hornung1997homogenization}
Ulrich Hornung.
\newblock {\em Homogenization and porous media}, volume~6.
\newblock Springer Science \& Business Media, 1997.

\bibitem{hornung2012homogenization}
Ulrich Hornung.
\newblock {\em Homogenization and porous media}, volume~6.
\newblock Springer Science \& Business Media, 2012.

\bibitem{hou1997multiscale}
Thomas~Y Hou and Xiao-Hui Wu.
\newblock A multiscale finite element method for elliptic problems in composite
  materials and porous media.
\newblock {\em Journal of computational physics}, 134(1):169--189, 1997.

\bibitem{le2014msfem}
Claude Le~Bris, Fr{\'e}d{\'e}ric Legoll, and Alexei Lozinski.
\newblock An msfem type approach for perforated domains.
\newblock {\em Multiscale Modeling \& Simulation}, 12(3):1046--1077, 2014.

\bibitem{lu2021uniform}
Yong Lu.
\newblock Uniform estimates for stokes equations in a domain with a small hole
  and applications in homogenization problems.
\newblock {\em Calculus of Variations and Partial Differential Equations},
  60:1--31, 2021.

\bibitem{muljadi2015nonconforming}
Bagus~Putra Muljadi, Jacek Narski, Alexei Lozinski, and Pierre Degond.
\newblock Nonconforming multiscale finite element method for stokes flows in
  heterogeneous media. part i: methodologies and numerical experiments.
\newblock {\em Multiscale Modeling \& Simulation}, 13(4):1146--1172, 2015.

\bibitem{weinan2003heterogenous}
E~Weinan and Bjorn Engquist.
\newblock The heterognous multiscale methods.
\newblock {\em Communications in Mathematical Sciences}, 1(1):87--132, 2003.

\bibitem{wolf2022homogenization}
Sylvain Wolf.
\newblock Homogenization of the stokes system in a non-periodically perforated
  domain.
\newblock {\em Multiscale Modeling \& Simulation}, 20(1):72--106, 2022.

\bibitem{xie2024cempoisson}
Wei Xie, Yin Yang, Eric Chung, and Yunqing Huang.
\newblock Cem-gmsfem for poisson equations in heterogeneous perforated domains.

\bibitem{yosifian1997some}
GA~Yosifian.
\newblock On some homogenization problems in perforated domains with nonlinear
  boundary conditions.
\newblock {\em Applicable Analysis}, 65(3-4):257--288, 1997.

\end{thebibliography}

\end{document}